\documentclass[graybox]{svmult}

\usepackage{type1cm}        
\usepackage{makeidx}         
\usepackage{graphicx}        
\usepackage{multicol}        
\usepackage[bottom]{footmisc}

\usepackage{newtxtext}       %
\usepackage[varvw]{newtxmath}       
\usepackage{physics} 

\makeindex             

\newcommand{\matid}{\ensuremath{\mathbf{{I}}}}
\newcommand{\bG}{\ensuremath{\mathbf{G}}}
\newcommand{\bx}{\ensuremath{\mathbf{x}}}
\newcommand{\bb}{\ensuremath{\mathbf{b}}}
\newcommand{\bz}{\ensuremath{\mathbf{z}}}
\newcommand{\br}{\ensuremath{\mathbf{r}}}
\newcommand{\bA}{\ensuremath{\mathbf{A}}}
\newcommand{\bH}{\ensuremath{\mathbf{H}}}
\newcommand{\by}{\ensuremath{\mathbf{y}}}
\newcommand{\bL}{\ensuremath{\mathbf{L}}}
\newcommand{\spann}{\ensuremath{\operatorname{span}}}
\newcommand{\bu}{\ensuremath{\mathbf{u}}}
\newcommand{\bW}{\ensuremath{\mathbf{W}}}
\newcommand{\bv}{\ensuremath{\mathbf{v}}}
\newcommand{\bP}{\ensuremath{\mathbf{P}}}
\newcommand{\bw}{\ensuremath{\mathbf{w}}}
\newcommand{\bX}{\ensuremath{\mathbf{X}}}
\newcommand{\bY}{\ensuremath{\mathbf{Y}}}

\begin{document}

\title*{Do \textit{you} precondition on the left or on the right?} 
\author{Nicole Spillane\orcidID{0000-0001-5067-946X} and\\ Pierre Matalon\orcidID{0000-0002-4288-2289} and \\ Daniel B.~Szyld\orcidID{0000-0001-8010-0391}}
\institute{Nicole Spillane \at CNRS, CMAP, \'Ecole polytechnique, Route de Saclay, 91128 Palaiseau cedex, France \email{nicole.spillane@cnrs.fr}
\and Pierre Matalon \at CNRS, CMAP, \'Ecole polytechnique, Route de Saclay, 91128 Palaiseau cedex, France \email{pierre.matalon@polytechnique.edu}
\and Daniel B.~Szyld \at Department of Mathematics, Temple University, Philadelphia, PA 19122, USA \email{szyld@temple.edu}}
%
%
\maketitle

\abstract*{
This work is a follow-up to a poster that was presented at the DD29 conference. Participants were asked the question `Do you precondition on the left or on the right?'. Here we report on the results of this social experiment. We also provide context on left, right and split preconditioning, share our literature review on the topic, and analyze some of the finer points. Two examples illustrate that convergence bounds can sometimes lead to misleading conclusions.}

\abstract{This work is a follow-up to a poster that was presented at the DD29 conference. Participants were asked the question: ``Do you precondition on the left or on the right?''. Here we report on the results of this social experiment. We also provide context on left, right and split preconditioning, share our literature review on the topic, and analyze some of the finer points. Two examples illustrate that convergence bounds can sometimes lead to misleading conclusions.}

\begin{figure}
\includegraphics[width=\textwidth]{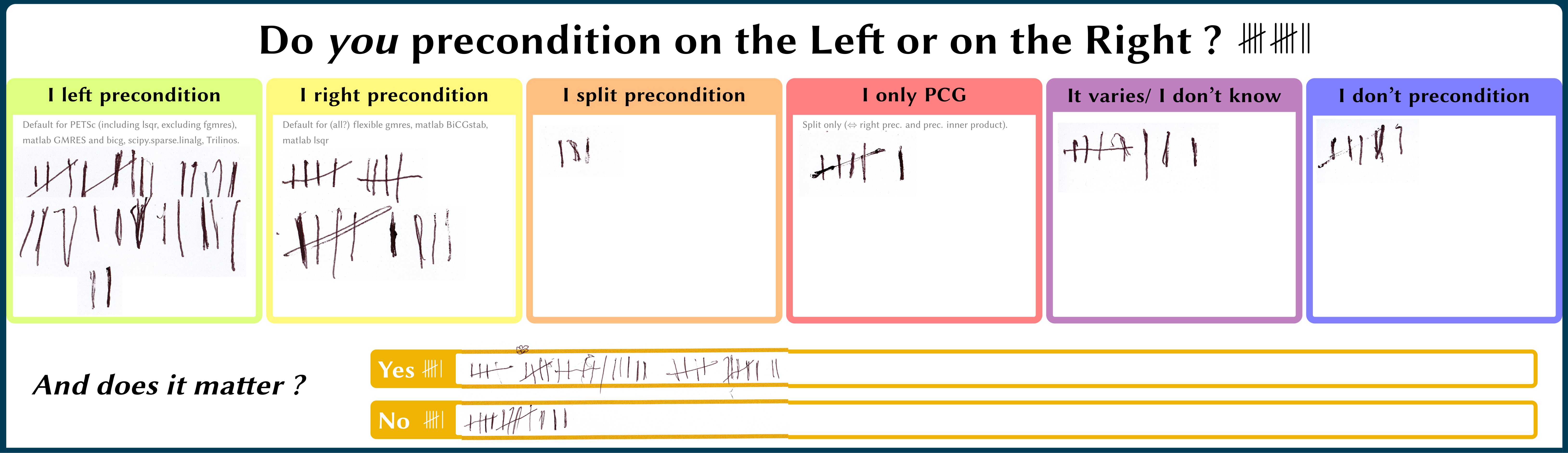}
\caption{Attendees of DD29 were asked to vote for their favourite way of preconditioning during the poster session. See Section~\ref{sec:DD29} and Table~\ref{tab:poll} for more details. }
\label{fig:poster}
\end{figure}

\section{Introduction}

This work considers Krylov subspace methods \cite[Chapters 6--9]{zbMATH01953444} for solving a linear system $\bA \bx = \bb$ given  $\bA \in \mathbb C^{n\times n}$ non-singular and  $\bb \in \mathbb C^n$. More precisely, we focus on the subset of Krylov subspace methods that minimize the residual. These methods start with an initial vector $\bx_0$, and thus the initial residual is $\br_0 := \bb - \bA \bx_0$. Then, each iteration $k\geq 1$ is characterized by an approximate solution $\bx_k$ which is the unique element such that the residual $\br_k := \bb - \bA \bx_k$ satisfies
\begin{equation}
\label{eq:min-krylov}
\| \br_k \|   =  \min_{\by \in  \bA \mathcal K_k(\bA , \br_0 )} \| \br_0 - \by \|,  
\end{equation} 
where the Krylov subspace is $ \mathcal K_k(\bA , \br_0 ) := \spann \{ \br_0 , \bA \br_0, \dots , \bA^{k-1} \br_0 \}$.

A well-established idea is to accelerate convergence by preconditioning. Given two non-singular matrices $\bH_L$ and $\bH_R$, the preconditioned system is
\begin{equation}
\label{eq:precsystem}
\bH_L \bA \bH_R\bu = \bH_L \bb; \text{ and } \bx = \bH_R \bu. 
\end{equation}
We refer to this system as being split-preconditioned by $(\bH_L,\bH_R)$ in which $\bH_L$ is the left preconditioner and $\bH_R$ is the right preconditioner. If $\bH_L = \matid$, then the system is right-preconditioned by $\bH_R$. If $\bH_R = \matid$, then the system is left-preconditioned by $\bH_L$. Both systems $\bA \bx = \bb$ and \eqref{eq:precsystem} are equivalent. The objective is to select the preconditioners in such a way that the Krylov subspace method is more efficient on~\eqref{eq:precsystem} than on $\bA \bx = \bb$. Applying the characterization \eqref{eq:min-krylov} to a Krylov subspace method for solving \eqref{eq:precsystem}  gives, with $\br_k = \bb - \bA \bx_k$,
\begin{equation}
\label{eq:min-pkrylov}
\|\bH_L \br_k \|   =  \min_{\by \in  \bA \mathcal K_k(\bH \bA , \bH \br_0 )} \| \bH_L(\br_0 - \by) \|,   
\end{equation} 
where $\mathcal K_k(\bH \bA , \bH \br_0 ) := \spann \{ \bH \br_0 , \bH \bA \bH \br_0, \dots , (\bH \bA)^{k-1}\bH \br_0 \}$ 
in which the combined preconditioner is defined by $\bH := \bH_R \bH_L$.

It is interesting to note that:
\begin{itemize}
\item the minimization space depends only on the choice of the combined preconditioner $\bH = \bH_R \bH_L$, 
\item the minimized norm depends only on the left preconditioner $\bH_L$. 
\end{itemize}

So far, we have not specified the norm $\| \cdot \|$. 
In fact, $\| \cdot \|$ can represent any weighted norm $\| \by \|_\bW := \sqrt{\by^* \bW \by}$ with $\bW$ Hermitian positive definite. 
Using such a norm yields the so-called \emph{weighted} GMRES method.
Nonetheless, we can choose the Euclidean norm without loss of generality. 
This is justified by the fact that \emph{weighting is preconditioning by similarity transformations} \cite{GL}: unpreconditioned GMRES using the norm $\| \cdot \|_\bW$ with $\bW := \bL {\bL^*}$ is equivalent to GMRES using the Euclidean norm and split-preconditioned by $({\bL^*}, {\bL^{-*}})$.
Consequently, the paper can be formulated either from the perspective of preconditioning or of weighted norms.
Since the preconditioning viewpoint is more familiar to most GMRES practitioners, we adopt this view, together with the Euclidean norm. 

Before setting weighted GMRES aside, we want to emphasize that the triplet $(\bH_L$, $\bH_R$, $\bW)$ contains redundancies that lead to equivalent configurations. 
For example, \emph{(i)} GMRES left-preconditioned by $\bH$ in norm $\| \cdot \|_\bW$ is equivalent to GMRES right-preconditioned by $\bH$ in norm  $\| \cdot \|_{\bH^*\bW\bH}$; 
\emph{(ii)} if $\bH := \bL {\bL^*}$ is Hermitian positive definite, then GMRES right-preconditioned by $\bH$ in norm $\| \cdot \|_\bH$ is equivalent to GMRES left-preconditioned by $\bH$ in norm $\| \cdot \|_{{\bH}^{-1}}$, which is also equivalent to GMRES split-preconditioned by $({\bL^*}, \bL)$ in Euclidean norm.
In practice, these equivalences offer significant flexibility. In particular, \emph{(ii)} shows how split preconditioning by the Cholesky factors of a Hermitian positive definite preconditioner can be realized without explicitly forming those factors.
The preconditioned conjugate gradient method exploits this technique to preserve symmetry \cite[Section 9.2.1]{zbMATH01953444}.
Weighted GMRES is also valuable for the theoretical analysis of Krylov methods.
Firstly, it creates a unifying framework that includes the conjugate gradient (CG) method by remarking that, if $\bA$ is Hermitian positive definite, then CG is simply GMRES expressed in the $\bA^{-1}$-norm.
Furthermore, many analytical results are established in weighted norm (see, \textit{e.g.}, \cite{zbMATH00036024,zbMATH07860856,zbMATH07931226,arXiv:2504.05723}), leaving open whether they extend to standard Euclidean-norm GMRES.
Addressing this latter question is one of the central motivations of the present work. This question has not been fully answered yet.

In the remainder of this work we compare left, right and split preconditioning for unweighted GMRES, which means that $\| \cdot\|$ is the Euclidean norm. We set out to determine which is preferable in terms of convergence and quality of the solution. In order to answer the question, Section~\ref{sec:litt} presents a literature review from which we extract the key points. Section~\ref{sec:analysis} analyzes some of these in more detail. In a more unconventional fashion, Section~\ref{sec:DD29} presents the results of a social experiment conducted at the DD29 conference on Domain Decomposition Methods where the participants were asked about their preferred way of preconditioning. 

\section{Literature Review}
\label{sec:litt}

We have turned to the literature and found surprisingly few concrete answers. In what follows we summarize what we found in the following references:
\begin{itemize}
\item {\textit{Iterative methods for sparse linear systems} by Y.~Saad \cite{zbMATH01953444}}
\item {\textit{Preconditioning} by A. J. Wathen \cite{zbMATH06443652} }
\item {\textit{Preconditioning and convergence in the right norm} by A. J. Wathen} \cite{zbMATH05198276} 
\item {\textit{Preconditioning techniques for large linear systems: A survey} by M. Benzi \cite{zbMATH01901695}} 
\item {\textit{Algorithms for sparse linear systems} by J. Scott and M.  Tůma  \cite{zbMATH07687320}}
\item {\textit{Any nonincreasing convergence curves are simultaneously possible for GMRES and weighted GMRES, as well as for left and right preconditioned GMRES} by P. Matalon and N. Spillane \cite{arXiv:2506.17193}}
\end{itemize}

\textbf{ \textit{Iterative methods for sparse linear systems} by Y.~Saad \cite{zbMATH01953444}}
\\
({Page 271})\footnote{The selected quotes also appear in the first edition of the book published in 1996 with a different page numbering. {In both cases the quote can be found in Section 9.3.}} \emph{
``A question arises on the differences between the right, left, and split preconditioning options. The fact that different versions of the residuals are available in each case may affect the stopping criterion and may cause the algorithm to stop either prematurely or with delay. This can be particularly damaging in case [$\bH$] is very ill-conditioned. The degree of symmetry, and therefore performance, can also be affected by the way in which the preconditioner is applied. For example, a split preconditioner may be much better if $\bA$ is nearly symmetric. Other than these two situations, there is little difference generally between the three options.''
}
 
({Page 271}) \emph{``When comparing the left, right, and split preconditioning options, a first observation to make is that the spectra of the three associated operators [$\bH\bA$, $\bA \bH$, and $\bH_L \bA \bH_R $]
are identical. Therefore, in principle one should expect convergence to be similar, although, as is known, eigenvalues do not always govern convergence.''
}

({Page 272})\emph{``In most practical situations, the difference in the convergence behavior of the two approaches is not significant. The only exception is when [$\bH$] is ill-conditioned which could lead to substantial differences.''
}

\medskip
\textbf{ \textit{Preconditioning} by A. J. Wathen \cite{zbMATH06443652} }
\\
(Pages 356-357) \emph{``One aspect which arises in the non-symmetric but not the symmetric case is the possibility to left-precondition [...], to right-precondition [...], or, if [$\bH$] is available in split form [$\bH = \bH_R \bH_L$], to use [split-preconditioning]. In every case the coefficient matrix is generally again non-symmetric, and obvious similarity transformations involving [$\bH$ or $(\bH_L, \bH_R)$] show that all three coefficient matrices are mathematically similar and so have the same eigenvalues. There is little evidence, however, that any one of these is better than the others, in terms of its effect on the convergence rate of a Krylov subspace method, even though, theoretically, this could be the case. Traditionally, right-preconditioning is favoured since then GMRES, for example, minimizes the residual $\br_k = \bb - \bA \bx_k$  for the original linear system; left-preconditioning would lead to minimization of the preconditioned residual [$\bH \br_k$]. It is certainly not always clear which is preferable: a good preconditioner fortunately seems to give fast convergence for all three forms in practice.''
}

\medskip
\textbf{\textit{Preconditioning and convergence in the right norm} by A. J. Wathen \cite{zbMATH05198276} }
\\
(Page 5) \emph{``[W]e show how preconditioning affects this balance between approximation and iteration error. In particular for widely used iterative methods which minimise the residual [...] we show that care is needed to avoid selection of preconditioners which apparently give rapid convergence, but which in fact merely distort the relevant norm so that poor solutions are achieved for all but extremely small convergence tolerances.''}

\medskip
\textbf{\textit{Preconditioning techniques for large linear systems: A survey} by M. Benzi \cite{zbMATH01901695}}
\\
(page 420) \emph{
``Which type of preconditioning to use depends
on the choice of the iterative method, problem characteristics, and so forth. For example, with residual minimizing methods, like GMRES, right preconditioning is often used. In exact arithmetic, the residuals for the right-preconditioned system are identical to the true residuals $\br_k = \bb- \bA \bx_k$.
\\ 
Notice that the matrices [$\bH \bA$, $\bA \bH$ , and $\bH_L \bA \bH_R$]  are all similar and therefore have the same eigenvalues. If $\bA$ and [$\bH$] are spd, the convergence of the CG method will be the same (except possibly for round-off effects) in all cases. On the other hand, in the nonnormal case, solvers like GMRES can behave very differently depending on whether a given preconditioner is applied on the left or on the right (see the discussion in \cite[{page 272}]{zbMATH01953444}\footnote{One of the quotations above by Y. Saad with page number edited for correct edition of book}; see also \cite[page 66]{SLee} for a striking example).''}
\medskip

\textbf{\textit{Algorithms for sparse linear systems} by J. Scott an M. T\accent23 uma \cite{zbMATH07687320}}
(page 167) \emph{``The following result states that it is not possible to determine a priori which variant is the best.
\\
{\underline{Theorem} ([Mendelsohn \cite{zbMATH03124252}]: }
 Let  $\delta$ and $\Delta$ be positive numbers. Then, for any $n \geq 3$, there exist non-singular $n \times n$ matrices $\bA$ and $\bH$ such that all the entries of $(\bH \bA - \matid)$ have absolute value less than $\delta$ and all the entries of $(\bA\bH - \matid)$ have absolute values greater than $\Delta$.
 \smallskip
\\ Nevertheless, the choice between left and right preconditioning is still important and may be based on the properties of the coupling of the preconditioner with the iterative method or on the distribution of the eigenvalues of A. The computed quantities that are readily available during a preconditioned iterative method depend on how the preconditioner is applied and this may influence the choice. These quantities may be used, for example, to decide when to terminate the iterations. An obvious advantage of right preconditioning is that in exact arithmetic, the residuals for the right preconditioned system are identical to the true residuals, enabling convergence to be monitored accurately. In some cases, the numerical properties of an implementation and/or the computer architecture may also play a part.''}
\medskip

\textbf{\textit{Any nonincreasing convergence curves are simultaneously possible for 
[...] left and right preconditioned GMRES} by P. Matalon and N. Spillane \cite{arXiv:2506.17193}}

(page 3) \emph{``
Consider two prescribed convergence curves of matching length:
\begin{itemize}
\item $\mathtt{r}_0 > \mathtt{r}_1 > \mathtt{r}_2 >  \dots $ for right preconditioned GMRES,
\item $\widetilde{\mathtt{r}_0} > \widetilde{\mathtt{r}_1} > \widetilde{\mathtt{r}_2} >  \dots $  for left preconditioned GMRES.
\end{itemize}
The first result is that there exists a system $A x = b$ and a preconditioner $H$ such that both convergence curves are realized. Additionally, the eigenvalues of $AH$ can be prescribed. 
[...] 
Besides highlighting that GMRES convergence does not solely depend on the eigenvalues, these  results show that for some cases, the decision to apply preconditioning on the left or on the right, or to apply weighted GMRES may lead to significant differences in convergence speed.
''}

\medskip 

From these quotes we draw the following conclusions:
\begin{enumerate}
\item The core difference between left and right preconditioning is the residual that gets minimized. This can be seen in the characterization \eqref{eq:min-pkrylov} (same $\bH$ but different~$\bH_L$).
\item
As a consequence of the previous observation, with right preconditioning, the residual norm can be monitored for free throughout the iterations. With left preconditioning, monitoring the preconditioned residual norm can lead to premature or delayed stopping of the algorithm in terms of the unpreconditioned residual norm. This important question, central to \cite{zbMATH05198276}, is also addressed, \textit{e.g.} in \cite{zbMATH01618501,zbMATH05626642}.
\item The left, right and split preconditioned matrices have the same spectrum so long as $\bH_R \bH_L$ is unchanged.
\item It is expected that convergence be faster for the matrix that is close to normality, or to symmetry. Ease of implementation could play a role too.
\item In practice, the behaviour is often not modified by switching the preconditioner between left and right. In fact there are few illustrations of opposite behaviours.
\item There is no easy answer. Indeed, any two convergence curves can be simultaneously prescribed for left and right preconditioned GMRES. And it is possible to find $\bA$ and $\bH$ such that the entries in $\bA \bH - \matid$ are arbitrarily large and those of $\bH \bA - \matid$ are arbitrarily small.
\end{enumerate}

\section{Analysis}
\label{sec:analysis}

Our own findings are in agreement with the literature review. It is theoretically possible for the behaviours of left and right preconditioned GMRES to be very different but usually nothing much happens. Below we detail three points. Before we go into those details we point out the relation between the eigenvectors in the three cases considered.
Let $\bv$ be an eigenvector of $\bH\bA$, then both $\bA \bv$ and $\bH^{-1} \bv$ are eigenvectors of $\bA\bH$ while both $\bH_L \bA \bv$ and  $\bH_R^{-1} \bv$ are eigenvectors of $\bH_L\bA\bH_R$ (all corresponding eigenvalues being equal).

\subsection{Left and right preconditioned GMRES can only differ if $\bH$ is ill-conditioned}

If right preconditioned GMRES produces residuals $\br_k$ and left preconditioned GMRES produces residuals $\widetilde{\br}_k$, they satisfy
\[
\|\br_k\| \leq \| \widetilde{\br}_k\| \leq \kappa(\bH) \| \br_k \|,
\]
where $\kappa(\bH) =  \sigma_{max}(\bH)/\sigma_{min}(\bH)$ is the condition number of $\bH$.  In other words, the residuals produced by left and right preconditioned GMRES, when measured in the Euclidean norm, cannot differ by a factor more than $\kappa(\bH)$. The proof is simply that
\[
\sigma_{min}(\bH)\|\br_k\| {\leq} \sigma_{min}(\bH) \| \widetilde{\br}_k\| {\leq}  \|\bH \widetilde{\br}_k \| {\leq} \|\bH \br_k\|{\leq} \sigma_{max}(\bH) \| \br_k \|,
\]
where the first and third inequalities are the optimality condition \eqref{eq:min-pkrylov} for right and left preconditioned GMRES, while the second and fourth inequalities are spectral estimates.

We can generalize this analysis to two split preconditioners $(\bH_L, \bH_R)$ and $(\widetilde{\bH}_L, \widetilde{\bH}_R)$ that satisfy $\bH_R \bH_L = \widetilde{\bH}_R \widetilde{\bH}_L$. The two minimization spaces are identical. Their corresponding residuals $\br_k$ and $\widetilde{\br}_k$ satisfy
\[
\sigma_{min}( \bG )\|\bH_L \br_k\| {\leq} \sigma_{min}( \bG) \|\bH_L \widetilde{\br}_k\| {\leq}  \|\widetilde{\bH}_L \widetilde{\br}_k \| {\leq} \|\widetilde{\bH}_L \br_k\| {\leq} \sigma_{max}( \bG   ) \|\bH_L \br_k \|,
\]
with $\bG = \widetilde{\bH}_L \bH_L^{-1} $. 
Again, the first and third inequalities are the optimality conditions, while the second and fourth inequalities are spectral estimates. It follows that: 
\[
\|\bH_L \br_k\| \leq \|\bH_L \widetilde{\br}_k\| \leq \kappa( \widetilde{\bH}_L \bH_L^{-1}) \| \bH_L \br_k \|.
\]

Applied to right preconditioning by $\bH$ (by setting $(\bH_L, \bH_R) = (\matid , \bH)$) versus split preconditioning by $(\widetilde{\bH}_L, \bH {\widetilde{\bH}_L}^{-1})$ we obtain 
\[
\| \br_k\| \leq \| \widetilde{\br}_k\| \leq \kappa( \widetilde{\bH}_L ) \| \br_k \|,
\]
so that right and split preconditioning can differ in Euclidean residual norm only by a factor $\kappa( \widetilde{\bH}_L )$. 

\subsection{Different convergence bounds, same behaviour (1/2)}

The striking example cited by M.~Benzi in \cite{zbMATH01901695} is from the PhD dissertation of S.~Lee \cite[page 66]{SLee}, and is adapted from the work of 
 A.~S.~Householder \cite[page 96]{zbMATH03516529}. It consists in a matrix $\bA$ and a preconditioner $\bH$ such that $\bH$ is a poor right inverse and a good left inverse in two senses: entrywise and field of values. 

The example is as follows. Let $\bA$ be a real matrix with eigenvalues in a real positive interval $[\mu, \lambda]$ with $0 < \mu < \lambda$. Let $\bu$ (respectively $\bv$) be a unit right (respectively left) eigenvector associated with $\lambda$ (respectively $\mu$),  
\textit{i.e.},
\[
\bA \bu = \lambda \bu; \quad \bv^\top \bA = \mu \bv^\top; \quad \| \bu \| = \|\bv \| = 1; \text{ and } \mu \neq \lambda.
\]
The proposed preconditioner is then a rank-one perturbation of $\bA^{-1}$ defined by 
\[
\bH :=  \bA^{-1} + \bu\bv^\top.
\]
It is easy to see that the preconditioned operators are 
$\bH\bA = (\bA^{-1} + \bu \bv^\top) \bA = \matid + \mu \bu \bv^\top$ and  $\bA\bH =\bA(\bA^{-1} + \bu \bv^\top)   =\matid + \lambda \bu \bv^\top$.
Moreover,
\[
\matid - \bH\bA = - \mu \bu \bv^\top \text{ and } \matid - \bA\bH = - \lambda \bu \bv^\top. 
\]
so the entries of $\matid - \bH\bA$ and $\matid - \bA\bH$ are in the ratio $\lambda/\mu$ which can be very large. This illustrates \cite[Theorem~9.2]{zbMATH07687320} (originally from \cite{zbMATH03124252}) cited in Section~\ref{sec:litt}. 

Moreover, the field of values of $\bH \bA$ is the circle centered at $1$ of radius $\mu/2$ and the field of values of $ \bA \bH$ is the circle centered at $1$ of radius $\lambda/2$ (see below for a proof). A well-known convergence bound for GMRES, the Elman bound \cite{ees1983},\cite[\S 5.3]{elman-phd} relates the rate of convergence of GMRES to the field of values of the operator. In particular, in the case where the field of values is a small disk around~$1$, convergence is provably fast. When $0$ is included in the field of values, Elman's bound does no longer apply. For this reason, if $\lambda > 1 \gg \mu$ it is tempting to conclude that $\bH$ is a poor right-hand side inverse and a good left-hand side inverse in terms of the field of values. 

However, when it comes to convergence of GMRES, the reality is quite different and this can be seen by observing that 
\[
(\bH \bA)^2 = (\matid + \mu \bu \bv^\top)^2 = \matid + 2  \mu \bu \bv^\top +  \mu^2 \bu \bv^\top \bu \bv^\top =   \matid + 2  \mu \bu \bv^\top = 2 \bH \bA - \matid,
\]
where the third equality is the result of $\bv \perp \bu$ since $\bu$ and $\bv$ are a pair of right and left eigenvectors associated with two different eigenvalues. Consequently, the grade of any vector with respect to $\bH\bA$ is 2 and, both left and right preconditioned GMRES converge in at most 2 iterations (for any right hand side and any starting vector). 

\textit{Proof of field of value formulae.} The field of values of $\bA\bH$ is the set of all complex numbers of the form
\[
\frac{\bz^* \bA \bH \bz }{ \bz^* \bz } =1 + \lambda \frac{\bz^* \bu \bv^\top \bz }{\bz^*\bz} = 1 + \lambda \frac{\bz^* \bu \cdot \bv^\top \bz }{\bz^*\bz}; \,  \bz \in \mathbb C^n,
\]
with $\bz^*$ the conjugate transpose of $\bz$.
Applying again that $\bv \perp \bu$, any $\bz \in \mathbb C^n$ can be split into $\bz = \xi_1 \bu + \xi_2 \bv + \bw$ with $\bw \perp \bu$, $\bw \perp \bv$, and $\xi_1, \, \xi_2 \in \mathbb C$. It follows that 
\[
\frac{\bz^* \bA \bH \bz }{ \bz^* \bz } =  1 + \lambda \frac{\overline{\xi_1} \xi_2}{|\xi_1|^2 + |\xi_2|^2 + \|\bw\|^2} .
\]
The inclusion of the circle in the field of values can be seen by setting $\bw = 0$, $\xi_1 = \mathrm{e}^{-\imag \theta}$ and $\xi_2 = \mathrm{e}^{\imag \theta}$  for all values of $\theta \in [0, \pi]$: 
\[
\frac{\bz^* \bA \bH \bz }{\bz^* \bz }= 1 + \frac{\lambda}{2} \mathrm{e}^{2 \imag \theta}.  
\]
The inclusion of the field of values in the circle is the result of 
\[
2|\overline{\xi_1} \xi_2 |\leq |\xi_1|^2 + |\xi_2|^2 \leq |\xi_1|^2 + |\xi_2|^2 + \|\bw\|^2 = \| \bz \|^2.  
\]
The proof for the field of values of $\bH\bA$ is almost identical.

\subsection{Different convergence bounds, same behaviour (2/2)}

The convergence bound \cite[Prop.\ 6.32]{zbMATH01953444} involving the eigenvalues and the conditioning of the eigenvectors is a well-known convergence bound for GMRES.
In this section, we build a family of matrices and preconditioners such that this bound is small when the preconditioner is applied on the right, large when it is applied on the left, and yet lead to the same convergence speed in practice.

Assume a non-singular, diagonalizable matrix $\bA := \bP \Lambda \bP^{-1}$, where $\Lambda$ is a diagonal matrix holding the eigenvalues of $\bA$, and $\bP$ is a matrix whose columns hold the eigenvectors of $\bA$.
The bound \cite[Prop.\ 6.32]{zbMATH01953444} for the GMRES method applied to the system $\bA\bx = \bb$ reads
\begin{equation} \label{eq:conv_bound}
    \frac{\|\br_k\|}{\|\br_0\|}  \leq  \epsilon^{(k)} \kappa(\bP),
\end{equation}
where 
\begin{equation*}
    \epsilon^{(k)} := \min_{q \in \mathbb{P}^k, q(0)=1} \max_{\lambda \in \sigma(\bA)} \abs{q(\lambda)}.
\end{equation*}
In $\epsilon^{(k)}$, $\mathbb{P}^k$ denotes the space of polynomial of degree at most $k$.
This term accounts for the spectrum of $\bA$, while the second term of the bound, $\kappa(\bP)$, accounts for the eigenvectors of $\bA$.
Note that $\kappa(\bP)$ is often seen as a measure of $\bA$'s departure from normality, in the sense that if $\bA$ is normal, then $\bP$ is unitary and $\kappa(\bP) = 1$.
On the other hand, if $\kappa(\bP)$ is large, it means that some eigenvectors are close to being linearly dependent, thus making $\bP^*\bP$ far from identity, and therefore $\bA^*\bA$ far from~$\bA\bA^*$.

Let us build a test case where left- and right-preconditioning lead to significantly different bounds.
Let $\Lambda$ be an invertible, well-conditioned, diagonal matrix.
Let $\bX$ be an ill-conditioned invertible matrix.
Let $\bY$ be a unitary matrix.
Now, define the matrix $\bA := \bX \Lambda \bY^{-1}$ and the preconditioner $\bH := \bY\bX^{-1}$.
We then have the left- and right-preconditioned matrices
$\bH\bA = \bY \Lambda \bY^{-1}$ and $\bA\bH = \bX \Lambda \bX^{-1}$.

Solving $\bA\bx=\bb$ with left-preconditioned GMRES applies GMRES to the system $\bH\bA\bx = \bH\bb$. Denoting by $\widetilde{\br}_k := \bH(\bb-\bA\bx_k)$ the $k^{th}$ residual produced, \eqref{eq:conv_bound} writes
\begin{equation*}
    \frac{\|\widetilde{\br}_k\|}{\|\widetilde{\br}_0\|}  \leq  \epsilon^{(k)}.
\end{equation*}
Indeed, since $\bY$ is unitary, $\kappa(\bY) = 1$.
We expect fast convergence in this case due to the sole dependency of the eigenvalues, which we have selected to be clustered.
Similarly, right-preconditioned GMRES yields the bound
\begin{equation*}
    \frac{\|\br_k\|}{\|\br_0\|}  \leq  \epsilon^{(k)} \kappa(\bX).
\end{equation*}
Since $\bX$ is ill-conditioned, this bound is large: we (naively) expect slow convergence.
However, we present below an example for which both cases converge almost at the same speed, with only small discrepancies in the successive residual norms produced.

\begin{figure}
    \begin{center}
    \includegraphics[width=0.7\textwidth]{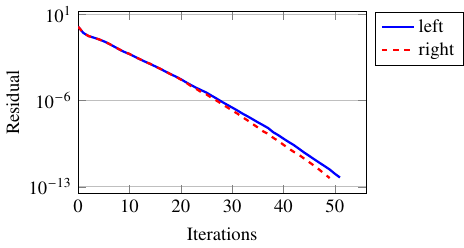}
    \caption{Example of left-right preconditioning leading to similar convergence behaviour while presenting significantly different convergence bounds.}
    \label{fig:example}
    \end{center}
\end{figure}

Experimental setup: $\bY$ is obtained as the orthogonal matrix produced by the QR factorization of a randomly generated square matrix.
$\bX$ is, first, initialized the same way. Then, denoting by $\bx_i$ the $i^{th}$ column of $\bX$, we replace $\bx_1$ by $\frac{1}{K}\bx_1 + \bx_n$, where $K$ is a large real number.
This way, $\bx_1$ is close to collinear to $\bx_n$, while preserving its strict linear independence from all other columns.
The coefficient $K$ drives the condition number of $\bX$.
$\bA$ and $\bH$ are defined accordingly, and $\bb$ is randomly generated.
Figure~\ref{fig:example} shows such an example with a system of order 100 and $K := 10^8$ ($\approx \kappa(\bX)$).
For right-preconditioning, we plot the residual norms $\|\br_k\|/\|\br_0\|$, and for left-preconditioning, we plot the preconditioned residual norms $\|\widetilde{\br}_k\|/\|\widetilde{\br}_0\|$. 
Note that the unpreconditioned residuals in this latter case sensibly produce the same convergence curve as the preconditioned residuals.

\section{A social experiment at DD29}
\label{sec:DD29}

\begin{table}
 \begin{center} 
\includegraphics[width=\textwidth]{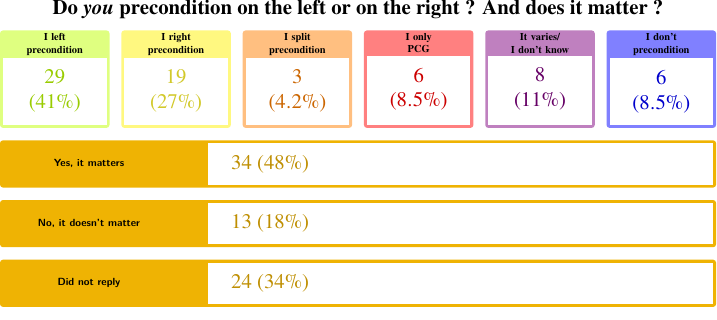}
\end{center}
\caption{Results of the poll. There are {71} replies to the first question. All percentages are with respect to these {71} participants.  Figure~\ref{fig:poster} shows how the questions were displayed on the poster.}
\label{tab:poll}
\end{table}

During the poster session at the DD29 conference on Domain Decomposition methods, the authors presented a poster which included a poll. The participants at the conference were asked to vote in favour of their favourite way of preconditioning. The possible answers and results of the poll are shown in Table~\ref{tab:poll}. There were 299 attendees at the conference according to the organizers. Of these, {$71$} people replied to the first question in the poll which is ``Do you precondition on the left or on the right''. The participation rate is  {$24\%$}. Significantly fewer people ({47}) replied to the second question which is ``And does it matter ?''. In the presentation of the results in Table~\ref{tab:poll} we included a ``Did not reply'' field to account for this.  

Preconditioned conjugate gradients (PCG) is given as a separate possible response to the first question since there is one overwhelmingly natural way of preconditioning the conjugate gradient method that is neither left- nor right-preconditioning \cite[Section 9.2.1]{zbMATH01953444} but preserves the short recurrence property.

It turns out that {$41\%$} of respondents apply left preconditioning. One possible explanation is that left preconditioning is the default choice in several libraries including PETSc (with the exception of fgmres but including lsqr), MATLAB's gmres and bicg, Python's scipy.sparse.linalg and Trilinos. Note that right preconditioning is the default for MATLAB solvers bicgstab and lsqr, and the necessary choice for all implementations of flexible GMRES. {In second place, 27\% of respondents apply right preconditioning.} The argument in favour of right preconditioning is to minimize and monitor the residual in the Euclidean norm\footnote{Even though the poll was not conducted in any other scientific meetings, the conversations that the first author has had with many scientists led her to believe that right preconditioning would have been the winner at a numerical linear algebra conference.}.  Finally, there are {2.6} times more people who consider the choice between left and right preconditioning to be important than those who don't. 

Both authors that were present at the conference responded that they left-precondition and that yes, it does matter. Our reasoning is that if $\bH$ is a good preconditioner for $\bA$, then the left preconditioned residual is a reasonable measure of the error. Since we are aware that the preconditioned residual norm may be much larger than the unpreconditioned residual norm, we do also compute the unpreconditioned residual (but not at every iteration) so that the solver does not stop prematurely.

\section{Conclusion}

The purpose of this work is to question whether there is a best way to position a given preconditioner. We have done our best to reflect existing considerations from the literature as well as to offer some enlightenments. Our examples show that convergence bounds that predict possible opposite behaviours between left and right preconditioned GMRES are not always correct predictions. If there is a norm that is particularly useful or meaningful for a given problem then it should be the norm that is minimized by GMRES. Otherwise, either choice is valid but GMRES users should be aware of what they are doing and perhaps compute several different residual norms for the output solution.

\begin{acknowledgement}
The work of the first two authors was funded in part by ANR project DARK (research grant ANR-24-CE46-1633).
\end{acknowledgement}

\end{document}